%% file: time_freezing_hysteresis_final.tex
\documentclass[journal,twoside,web]{ieeecolor}
\usepackage{lcsys}
\newcommand\copyrighttext{%
	\footnotesize \textcopyright 2022 IEEE. Personal use of this material is permitted. Permission from IEEE must be obtained for all other uses, in any current or future media, including reprinting/republishing this material for advertising or promotional purposes, creating new collective works, for resale or redistribution to servers or lists, or reuse of any copyrighted component of this work in other works.}
\newcommand\copyrightnotice{%
	\begin{tikzpicture}[remember picture,overlay]
		\node[anchor=south,yshift=5pt] at (current page.south) {\fbox{\parbox{\dimexpr\textwidth-\fboxsep-\fboxrule\relax}{\copyrighttext}}};
	\end{tikzpicture}%
}

\usepackage{cite}
\usepackage{amsmath,amssymb,mathtools,amsfonts}
\usepackage{graphicx}
\usepackage{textcomp}

\usepackage{hyperref}
\usepackage{multirow}
\usepackage{color}
\usepackage{psfrag}
\usepackage{subfigure}

\usepackage{empheq}
\usepackage{tabularx}
\usepackage{bm} 
\usepackage{url}
\usepackage{psfrag}
\usepackage{float}
\usepackage{comment}

\usepackage{alphalph}

\usepackage[short]{optidef}

\usepackage{amsthm}

\newtheorem{theorem}{Theorem}

\newtheorem{proposition}[theorem]{Proposition}

\newtheorem{definition}[theorem]{Definition}

\theoremstyle{definition}

\renewcommand{\QED}{\QEDopen}

\usepackage{tikz}
\usepackage{pgfplots}
\usetikzlibrary{arrows.meta}
\newlength\fwidth
\newlength\fheight
\usetikzlibrary{arrows}

\newcommand{\R}{{\mathbb{R}}}

\newcommand{\ts}{t_{\mathrm{s}}}
\newcommand{\tf}{t_{\mathrm{f}}}
\newcommand{\taus}{\tau_{\mathrm{s}}}
\newcommand{\dd}{\mathrm{d}}
\newcommand{\taur}{\tau_{\mathrm{r}}}
\newcommand{\tauf}{\tau_{\mathrm{f}}}
\newcommand{\tauj}{\tau_{\mathrm{jump}}}

\newcommand{\A}{{\mathrm{A}}}
\newcommand{\B}{{\mathrm{B}}}

\newcommand{\wone}{0}
\newcommand{\wtwo}{1}
\newcommand{\Tjmp}{\mathcal{T}_{\mathrm{jump}}}

\begin{document}
\title{Continuous Optimization for Control of Hybrid Systems with Hysteresis via Time-Freezing}
\author{Armin Nurkanovi\'c and Moritz Diehl
\thanks{This research was supported by the DFG via Research Unit FOR 2401 and project 424107692 and by the EU via ELO-X 953348.}
	\thanks{Armin Nurkanovi\'c is with the Department of Microsystems Engineering (IMTEK), University of Freiburg, Germany, Moritz Diehl is with the Department of Microsystems Engineering (IMTEK) and Department of Mathematics, University of Freiburg, Germany,  \texttt{\{armin.nurkanovic,moritz.diehl\}@imtek.uni-freiburg.de}}
}
\maketitle
\copyrightnotice
\thispagestyle{empty} 
\begin{abstract}
This article regards numerical optimal control of a class of hybrid systems with hysteresis using solely techniques from nonlinear optimization, without any integer variables.
Hysteresis is a rate independent memory effect which often results in severe nonsmoothness in the dynamics.
These systems are not simply Piecewise Smooth Systems (PSS); they are a more complicated form of hybrid systems. We introduce a time-freezing reformulation which transforms these systems into a PSS.
From the theoretical side, this reformulation opens the door to study systems with hysteresis via the rich tools developed for Filippov systems. 
From the practical side, it enables the use of the recently developed Finite Elements with Switch Detection~\cite{Nurkanovic2022}, which makes high accuracy numerical optimal control of hybrid systems with hysteresis possible.
We provide a time optimal control problem example and compare our approach to mixed-integer formulations from the literature.
 \end{abstract}
\vspace{-0.3cm}

\section{Introduction}\label{sec:introduction}
\IEEEPARstart{H}{ysteresis} occurs in many physical systems, e.g., ferromagnetism, plasticity, superconductivity, phase transitions, but also in feedback control, e.g., thermostats \cite{Visintin1994,Lunze2009}.
Hysteresis effects in dynamic systems are modeled with nonsmooth differential equations. This paper focuses on transforming some classes of systems with hysteresis into piecewise smooth system (PSS) and numerically solving optimal control problems (OCP) with PSS. We leverage recent advances in numerical optimal control of PSS, namely we use the FESD method~\cite{Nurkanovic2022}. 

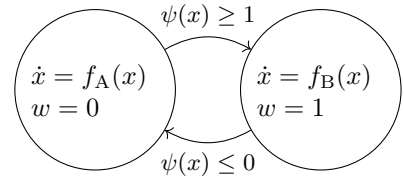
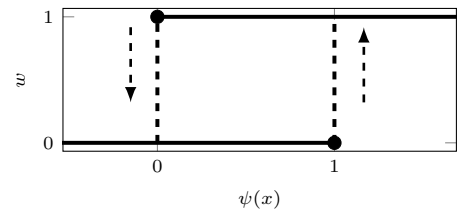
\begin{figure}
	\vspace{-0.05cm}					
	\centering
	\subfigure[center][Finite automaton of a hybrid system with hysteresis]{\input{img/finite_automaton2.tikz}\label{fig:finite_automation}}
	\hfill
	\subfigure[center][Example hysteresis characteristic $(w,\psi(x))$]{\input{img/hysteresis_char.tikz}
	\label{fig:hysteresis_char}}
	\caption{Hybrid system with hysteresis.}
	\label{fig:hysteresis_system}
	\vspace{-0.50cm}					
\end{figure}

A hybrid system with hysteresis can be represented as a \textit{finite automaton} \cite{Lunze2009} which has two modes of operation described by $f_{\mathrm{A}}(x)$ and $f_{\mathrm{B}}(x)$, cf. Fig. \ref{fig:finite_automation}.
If the system operates in mode $\A$ with $\dot{x}=f_{\mathrm{A}}(x)$ and if $\psi(x)\geq 1$, it switches to mode $\mathrm{B}$ with $\dot{x}=f_{\B}(x)$. On the other hand, if it operates in mode $\B$ and if $\psi(x)\leq 0$, it switches to mode $\A$. This is a typical hysteresis behavior given by the characteristic in Fig. \ref{fig:hysteresis_char}, which is often called the \textit{delayed relay operator} \cite{Brogliato2020}.  
The dynamics of the system depend on the value of $w(t)$ and the scalar \textit{switching function} $\psi(x)$. Notably, for $\psi(x)\in [0,1]$ the function $w(t)$ can be 0 or 1.
 \color{black}

There are several other related characteristics, e.g. the dashed lines in Fig. \ref{fig:hysteresis_char} could be solid, or the resulting polygon in the middle of the plot might be tilted. In all these cases the characteristic can be readily represented via a linear complementarity problem \cite{Vandenberghe1989} and the nonsmooth dynamic system recast into a Dynamic Complementarity System (DCS). However, it is an open question if this DCS is a PSS.

{Control of systems with hysteresis relying on Filippov solutions was studied in, e.g., \cite{Ceragioli2011,Bagagiolo2019}.} In control theory, systems with hysteresis are often studied via the hybrid systems framework which uses integer state and control variables \cite{Lunze2009,Bemporad1999b,Avraam2000}. Hence, in an optimal control context this requires solving  Mixed Integer Optimization Problems (MIOP). They can be solved efficiently in case of discrete time linear hybrid systems \cite{Bemporad1999b} where MILP or MIQP formulations can be found. However, as soon as the junction times need to be determined precisely or non-linearity is present, e.g., in time optimal control problems, solving MIOP can become arbitrarily difficult. 
On the other hand, the nonsmoothness can be modeled with complementarity constraints \cite{Nurkanovic2020} and one must solve only nonsmooth Nonlinear Programs (NLP). However, standard time-stepping schemes for DCS have only first order accuracy and result necessarily in wrong numerical sensitivities and artificial local minima \cite{Stewart2010,Nurkanovic2020}. 

The \textit{time-freezing} reformulation transforms systems with state jumps into PSS and was first introduced in \cite{Nurkanovic2021,Nurkanovic2021c}. This paper introduces a {time-freezing} reformulation to transform systems represented with the {finite automaton} in Fig. \ref{fig:finite_automation} into PSS. 
Here, the main idea is to regard $w(t)$ as a continuous differential state. However, $w(t)$ exhibits jump discontinuities in time at $(0,1)$ and $(1,0)$, which can be interpreted as a \textit{state jump law}. 
As in \cite{Nurkanovic2021,Nurkanovic2021c}, we introduce \textit{auxiliary dynamic systems} and a \textit{clock state}. The auxiliary dynamic systems evolve in regions which are prohibited for the initial system and their trajectory endpoints satisfy the state jump law. Additionally, the evolution of the clock state is frozen during the evolution of the auxiliary systems. By regarding only the parts of $w(\cdot)$ when the clock state was evolving, we recover the original discontinuous solution. Note that the resulting time-freezing system is now a PSS, since the only remaining jump discontinuities are in the system's dynamics but not in the state anymore. For high accuracy numerical optimal control of PSS we use the FESD method \cite{Nurkanovic2022}. An implementation is available in the open source software package \texttt{NOSNOC} \cite{Nurkanovic2022b,Nurkanovic2022c}.

\paragraph*{Contribution} We present a time-freezing reformulation for a class of hybrid systems with hysteresis, which transforms them into PSS. 
Constructive ways for finding the auxiliary dynamics needed in time-freezing are provided. Solution equivalence between the initial hybrid and time-freezing PSS are proven. From the theoretical side, this contribution enables one to treat hybrid systems with hysteresis with the tools for PSS and Filippov systems \cite{Filippov2013}.
From the practical side, the highlight of this paper is that we can solve OCP with systems with hysteresis with high accuracy and without the use of any integer variables. The OCP discretized via FESD result in Mathematical Programs with Complementarity Constraints (MPCC). With appropriate reformulations the MPCC can often be solved by only a few NLP solves \cite{Anitescu2007}, i.e., the highly nonsmooth and nonlinear OCP are solved by purely derivative based algorithms. A time optimal control problem of a hybrid system with hysteresis and illustrates theoretical and algorithmic developments.  We compare the continuous optimization-based FESD method to mixed integer solution strategies.\color{black}
\paragraph*{Outline} Section \ref{sec:basic_defintions} gives some basic definitions on hybrid systems with hysteresis and PSS. In Section \ref{sec:time_freezing} we develop the time-freezing reformulation for a class of hybrid systems with hysteresis and provide a simple tutorial example. Section \ref{sec:solution_relationship} formalizes the relation between time-freezing PSS and hysteresis systems. Finally, Section \ref{sec:numerical_example} contains a numerical example and Section \ref{sec:conclusion} concludes the paper.
\paragraph*{Notation}
  
For the \textit{phyisical} time derivative of a function $x(t)$ we use $\dot{x}(t) \coloneqq \frac{\dd x}{\dd t}(t)$ and for the \textit{numerical} time derivative of $y(\tau)$ we use ${y}'(\tau)\coloneqq\frac{\dd y}{\dd\tau}(\tau)$. \color{black}
The matrix ${I_n \in \R^{n\times n}}$ is the identity matrix, and $\mathbf{0}_{m,n} \in \R^{m\times n}$ is the zero matrix.  The concatenation of two column vectors $a\in \R^m$, $b\in \R^n$ is denoted by $(a,b)\coloneqq [a^\top,b^\top]^\top$. The concatenation of several column vectors is defined in an analogous way. The closure of a set $C$ is denoted by $\overline{C}$, its boundary as $ \partial C$   and ${\textmd{conv}}(C)$ is its convex hull. \color{black}
\section{Basic Definitions: Hybrid Systems with Hysteresis and Filippov Systems} \label{sec:basic_defintions} 
In this section we provide some of the basic definitions and notation for PSS and hybrid systems. 
\subsection{PSS and Filippov Systems}\label{sec:pss_and_filippov} 
We regard PSS of the following form
\begin{align} \label{eq:pws1}
		\dot{x} & =f_i(x),\ \text{if} \; x \in R_i \subset \R^{n_x},\ i \in\! \mathcal{I} \coloneqq\! \{ 1,\ldots,m \},
\end{align}
with  regions $R_i\subset \R^{n_x}$ and associated dynamics $f_i(\cdot)$, which are smooth functions on an open neighborhood of $\overline{R}_i$ and $m$ is a positive integer. Note that in general the right hand side (r.h.s.)  of \eqref{eq:pws1} is discontinuous in $x$.
We assume that $R_i$ are disjoint, nonempty, connected and open sets. They have piecewise-smooth boundaries $\partial R_i$.
  Moreover, we assume that $\overline{\bigcup\limits_{i\in \mathcal{I}} R_i} = \R^{n_x}$ and that $\R^{n_x} \setminus \bigcup\limits_{i\in\mathcal{I}} R_i$ is a set of measure zero. \color{black} Note that the dynamics are not defined on $\partial R_i$ and to have a meaningful solution concept for the PSS \eqref{eq:pws1} we regard the Filippov convexification of it \cite{Filippov2013}. The ODE \eqref{eq:pws1} with a disconitous r.h.s. is replaced by a Differential Inclusion (DI) whose r.h.s. is a convex and bounded set. Due to the assumed structure of the sets $R_i$, if $\dot{x}$ exists, functions $\theta_i(\cdot)$ which serve as convex multipliers can be introduced and the Filipov DI for \eqref{eq:pws1} reads as \cite{Nurkanovic2022}
\begin{align}\label{eq:FilippovDI_with_multiplers}
	\begin{split}
		\dot{x}  \in    F_{\mathrm{F}}(x) &= \Big\{ \sum_{i\in \mathcal{I}}
		f_i(x) \, \theta_i  \mid \sum_{i\in \mathcal{I}}\theta_i = 1,
		\ \theta_i \geq 0,\\
		&0= \theta_i \  \mathrm{if} \;  x \notin \overline{R_i}, 
		\forall  i  \in \mathcal{I} \Big\}.
	\end{split}
\end{align}
Note that in the interior of the regions $R_i$ the \textit{Filippov set} $F_\mathrm{F}(x)$ is equal to $\{f_i(x)\}$ and on the boundary between regions we have a convex combination of the neighboring vector fields. The evolution of $x(\cdot)$ on region boundaries $\partial R_i$ are called \textit{sliding modes}. The sliding mode dynamics in Filippov's setting are implicitly defined by Differential Algebraic Equations (DAE) \cite{Filippov2013}.

%

\subsection{Hybrid Systems with Hysteresis} \label{sec:hysteresis_system}
We consider dynamic systems represented with the finite automaton in Fig. \ref{fig:finite_automation}:
 
\begin{align}\label{eq:ode_hysteresis}
	\dot{x} &= f(x,w) = \!(1-w)f_{\A}(x)+wf_{\B}(x),
\end{align}
where the $(w,\psi(x))$ characteristic is illustrated in Fig. \ref{fig:hysteresis_char}. 
\color{black}
For a uniformly continuous function $x(t)$ on $t \in [0,T]$ and a smooth $\psi(\cdot)$, there can be only finitely many oscillations  between $0$ and $1$. Consequently, the function $w(t)$ is piecewise constant and has only finitely many jumps between $0$ and $1$ \cite{Visintin1994}. 

The system in \eqref{eq:ode_hysteresis} has two modes of operation denoted by ${\A}$ and ${\B}$. In order to be able to simulate \eqref{eq:ode_hysteresis} for $t\in [0,T]$ with a given $x(0)  = x_0$ we must know $w(0)$ as well. This property is typical for systems with hysteresis. Furthermore, $w(\cdot)$ jumps between $\wone$ and $\wtwo$, hence we can describe it by an ODE with the state vector $z \coloneqq (x,w) \in \R^{n_x+1}$ which is associated with a state jump law. 
\begin{align}\label{eq:ode_hysteresis_with_w}
	\dot{z} &= (f(x,w),0),
\end{align}
accompanied by a state-jump law for $w(\cdot)$ at time-point $\ts$ which covers two scenarios:
\begin{enumerate}
	\item  if $w(\ts^{-}) = 0$ and $\psi(x(\ts^{-})) = 1$, then $x(\ts^+) = x(\ts^{-})$ and $w(\ts^{+}) = 1$,
	\item  if $w(\ts^{-}) = 1$ and $\psi(x(\ts^{-})) = 0$,  then $x(\ts^+) = x(\ts^{-})$ and $w(\ts^{+}) = 0$.
\end{enumerate}
Clearly, due to the state jump law the ODE \eqref{eq:ode_hysteresis_with_w} is not simply a PSS as \eqref{eq:pws1}. Throughout the paper we assume, given $x(0)$ and $w(0)$ that there exists a solution to the Initial Value Problem (IVP) associated with \eqref{eq:ode_hysteresis_with_w}.   A way to define a meaningful notion of solution for hybrid system as \eqref{eq:ode_hysteresis_with_w} is given in  e.g., \cite[Section 5.4]{Lunze2009} and sufficient conditions for well-posedness are provided  \cite[Theorem 5.4]{Lunze2009}.\color{black}

\section{The Time-Freezing Reformulation for Hybrid Systems with Hysteresis}\label{sec:time_freezing}
This section introduces the time-freezing reformulation for the system \eqref{eq:ode_hysteresis_with_w}. 
We define step-by-step the corresponding regions $R_i$ of the time-freezing PSS and give constructive ways to find vector fields associated to them. The section finishes with a tutorial example.
\subsection{The Time-Freezing System}
The main idea is to transform the state $w(t)$ which is a piecewise constant function of time into a continuous differential state on a different time domain. We call this new time domain the \textit{numerical time} and denote it by $\tau$. Instead of $t$ as in \eqref{eq:pws1}, $\tau$ will now be the time of the time-freezing PSS.
Moreover, we introduce a clock state $t(\tau)$ in the time-freezing PSS which we call \textit{physical time}. It grows whenever the systems evolves according to $f_{\A}(x)$ or $f_{\B}(x)$, i.e., $\frac{\dd t}{\dd \tau}(\tau) = 1$. Otherwise the physical time is frozen, i.e., $\frac{\dd t}{\dd \tau}(\tau) = 0$. In other words, the time is frozen whenever $w \notin \{0,1\}$. Consequently, the $w(\cdot)$ takes only discrete values in physical time, i.e., when $t(\tau)$ is evolving. 
 
The time-freezing PSS has the following state vector $y \coloneqq (x,w,t) \in \R^{n_y}, n_y = n_x+2$. In the sequel, we define its regions $R_i \subset \R^{n_y}$ and the associated vector fields $f_i(y)$. Some key observations can be made from Fig \ref{fig:hysteresis_char}. First, everything except the solid curve is prohibited for the system \eqref{eq:ode_hysteresis_with_w} in the $(\psi,w)-$ plane. We use this prohibited part of the state space to define auxiliary dynamics. Second, the evolution happens in a lower-dimensional subspace since $\dot{w} = 0$. This corresponds in Filippov's setting to {sliding modes}. Hence, we define the regions such that the evolution of the initial system \eqref{eq:ode_hysteresis_with_w} corresponds to sliding modes of the time-freezing PSS, i.e., it happens on region boundaries $\partial R_i$.

\begin{figure}[t]
	\centering
	\includegraphics[scale=0.82]{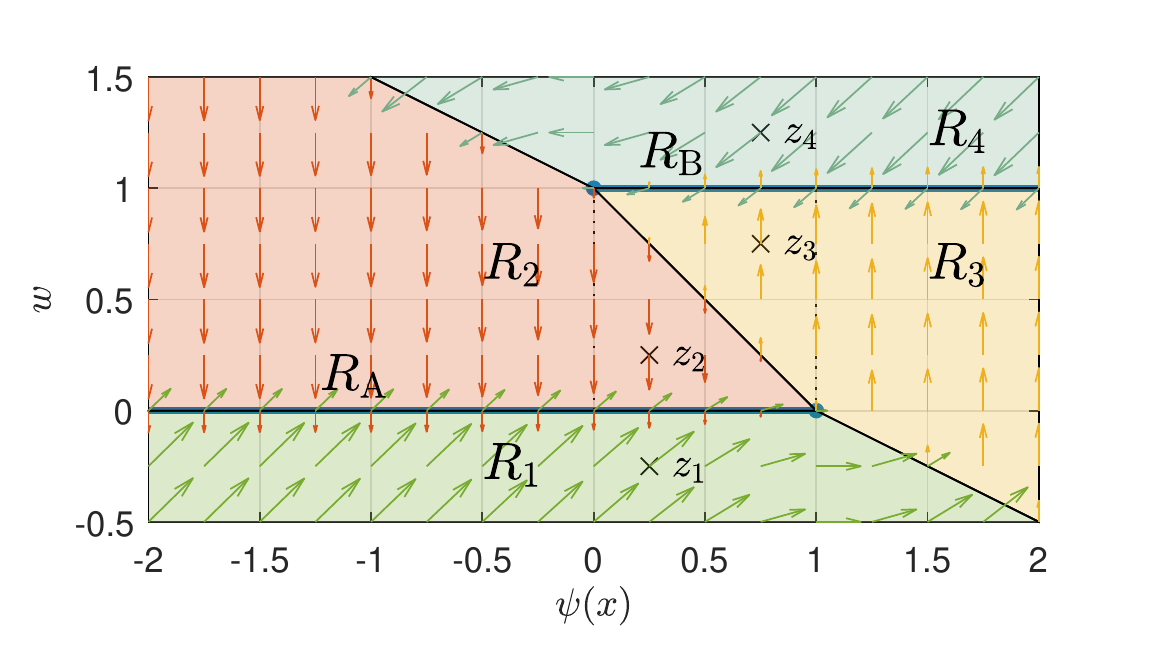}
	\vspace{-0.4cm}			
	\caption{Illustration of the partitioning of the state space in $(\psi(x),w)$-plane for the time-freezing PSS via Voronoi regions with the corresponding auxiliary and DAE-forming dynamic's vector fields. The Voronoi points $z_i,  i=1,\ldots,4$, are marked by the crosses.}
	\label{fig:vector_fields_hsyschar}			
	\vspace{-0.5cm}			
\end{figure}

A suitable partition of the $(\psi,w)-$ plane can be achieved with \textit{Voronoi regions}. The regions are defined as $R_i = \{ z \mid \|z-z_i \|^2 < \|z-z_j \|^2,\; j = {1,\ldots,4}, j\neq i \}$, $z = (\psi(x),w)$ with the points: $z_1 = (\frac{1}{4},-\frac{1}{4}), z_2 = (\frac{1}{4},\frac{1}{4}), z_3 = (\frac{3}{4},\frac{3}{4})$ and $z_4 = (\frac{3}{4},\frac{5}{4})$.
 An illustration of the regions is given in Fig. \ref{fig:vector_fields_hsyschar}, where the black solid lines denote the region boundaries.   This choice of $z_i$ defines regions such that their boundaries correspond to the feasible set of the original system \eqref{eq:ode_hysteresis_with_w}. Moreover, the space is split by the diagonal line between $R_2$ and $R_3$ such that we can define different auxiliary dynamics for the state jumps in both directions. One can make other choices for the points $z_i$ with the exact same properties. The proposed choice partitions the space symmetrically, cf. Fig. \ref{fig:vector_fields_hsyschar}. \color{black}
 The figure illustrates also the vector fields in the regions $R_i$ whose meaning is detailed below. It is important to note that the original system can only evolve at region boundaries $R_{\A}\coloneqq \{y \in \R^{n_y} \mid w =0,\ \psi(x)\leq1\} = \partial {R}_1 \cap \partial R_2$ and $R_{\B} \coloneqq \{y \in \R^{n_y} \mid w =1,\ \psi(x)\geq0\} =\partial R_3 \cap \partial R_4$. 

We exploit the interior of the regions $R_i$, $i = 1,\dots,4$ to define the needed auxiliary ODE. In what follows, in the regions $R_2$ and $R_3$ we define auxiliary dynamic systems whose trajectory endpoints satisfy the state jump law of \eqref{eq:ode_hysteresis_with_w}. In the regions $R_1$ and $R_4$ we will define so-called DAE-forming ODE \cite{Nurkanovic2021c}, which make sure that we obtain appropriate sliding modes on $R_{\A}$ and $R_{\B}$, which are described by index 2 DAE \cite{Filippov2013} and witch match the dynamics of the original system. The next definition formalizes the desired proprieties of an auxiliary ODE.

%
 
\begin{definition}[Auxiliary ODE]\label{def:auxiliary_ode}
The auxiliary ODE in regions $R_2$ and $R_3$ are denoted by $y' = f_{\mathrm{aux},\A}(y)$ and $y' = f_{\mathrm{aux},\B}(y)$, respectively. For every initial value $y(\taus) = y_{\mathrm{s}}$ such that $(w(\taus),\psi(x(\taus))=(1,0)$, for $y_{\mathrm{s}}\in {R}_{\B}$, (and $(w(\taus),\psi(x(\taus)) = (0,1)$ for $y_{\mathrm{s}}\in {R}_{\A}$, respectively) and for every well-defined and finite time interval $\Tjmp \coloneqq (\taus,\taur)$ with the length $\tauj \coloneqq \taur-\taus$, the auxiliary ODE satisfy the following properties: 
(i) $w(\tau) \in (0,1), \ \forall \tau \in \Tjmp$, 
(ii) $x(\taus)=x(\taur)$,
 and (iii) $w(\taur) = 0$ (or $w(\taur) = 1$). 
\end{definition}
 \color{black}
 In other words, we define an ODE whose trajectory endpoints on $\overline{\mathcal{T}}_{\mathrm{jump}}$ satisfy the state jump law associated with Eq. \eqref{eq:ode_hysteresis_with_w}, cf. Fig. \ref{fig:vector_fields_hsyschar}. The next proposition provides a constructive way to find an ODE with the above described properties.
\begin{proposition}[Auxiliary ODE]\label{prop:auxiliary_ode}
	Given an initial value $y(\taus) = y_{\mathrm{s}}$ such that $w(\taus) = 1$ and $\psi(x(\taus)) =0$, the ODE given by
	\begin{align}\label{eq:auxiliary_ode1}
		y'(\tau) &= f_{\mathrm{aux},\A} (y) \coloneqq (\textbf{0}_{n_x,1},-\gamma(\psi(x)-1),0),
	\end{align}
	where $\gamma: \R\to\R$ and $\gamma(x) = \frac{ax^2}{1+x^2}$ with $a>0$, is an auxiliary ODE defined in $R_2$.
	Similarly, for $y(\taus) = y_{\mathrm{s}}$ with $w(\taus) = 0$ and $\psi(x(\taus)) =1$, the ODE
	 
	\begin{align}\label{eq:auxiliary_ode2}
		y'(\tau) &= f_{\mathrm{aux},\B} (y) \coloneqq (\textbf{0}_{n_x,1}, \gamma(\psi(x)),0).
	\end{align}
	\color{black}
	is an auxiliary ODE in $R_3$. In both cases $\tauj =\frac{1}{\gamma(-1)}$.
\end{proposition}
\textit{Proof.} We prove the assertion for \eqref{eq:auxiliary_ode1}, since the second part follows similar lines. Since $x'(\tau) = \textbf{0}_{n_x,1}$ and $t'(\tau) = 0$ these two variables do not change their value, thus $\psi(x(\tau)) = \psi(x(\taus)) = 0$ and $t(\tau) = t(\taus)$ for $\tau \geq \taus$. Hence, we have $w'(\tau) = -\gamma(-1)<0$. By explicitly solving the ODE we obtain $w(\taur) = 0$ for $\taur = \taus + \tauj$, where $\tauj =\frac{1}{\gamma(-1)}$. All conditions of Definition \ref{def:auxiliary_ode}  are satisfied thus the proof is complete. 
\qed

We briefly discuss some of the proprieties of such an auxiliary ODE, since the are several ways to construct similar ODE. Loosely speaking, in Fig. \ref{fig:vector_fields_hsyschar} in $R_2$ the vector field should point in the negative $w$-detection and in $R_3$ in the positive $w$-direction, and be zero in all other directions.
Note that for $\psi(x) \in (0,1)$ the vector fields of the auxiliary ODE in both cases point away from the manifold defined $\mathcal{M} = \{y \in \R^{n_y} \mid w + \psi(x)-1 = 0\}$. In such scenarios, there is usually locally no unique solution to the associated Filippov DI, as the trajectory can leave $\mathcal{M}$ at any point in time \cite{Filippov2013}.   However, the system should never be initialized in this region, since this state is infeasible for the original system. \color{black} We show later that it can never reach this undesired state if initialized appropriately. Furthermore, the auxiliary ODE from Proposition \ref{prop:auxiliary_ode} have by construction the favorable property that they do not point away in both directions from $\mathcal{M}$ at the junction points $(0,1)$ and $(1,0)$. This is why the function $\gamma(\cdot)$ was introduced in the auxiliary ODE.
Another favorable property is, if the system is initialized with the wrong value for $w(\cdot)$ for $\psi(x)\notin (0,1)$ the auxiliary ODE will automatically reinitialize $w(\cdot)$ while the physical time is frozen, cf. Fig \ref{fig:vector_fields_hsyschar}.

We still need to define DAE-forming vector fields for the regions $R_1$ and $R_4$. 
These vector fields should be such that, together with the auxiliary dynamics in their respective regions, they results in sliding modes on $R_{\A}$ and $R_{\B}$ which match the dynamics of the initial system \eqref{eq:ode_hysteresis_with_w}. 
 
In a general PSS the vector fields are not defined on the region boundaries, thus we use Filippov's convexification \cite{Filippov2013} as defined in Eq. \eqref{eq:FilippovDI_with_multiplers}, and denote the Filippov set associated to the time-freezing PSS by $F_{\mathrm{TF}}(\cdot)$. The next proposition gives a constructive way to find the desired vector fields.
\color{black}
\begin{proposition}[DAE-forming ODE]\label{prop:dae_forming_ode}
Suppose the regions $R_2$ and $R_3$ are equipped with the vector fields $f_{\mathrm{aux},\A}(\cdot)$ and $f_{\mathrm{aux},\B}(\cdot)$ from Proposition \ref{prop:auxiliary_ode}, respectively. Let the region $R_1$ be equipped with the  ODE
	\begin{align}\label{eq:dae_forming_ode1}
		y' &= f_{\mathrm{DF},\A}(y) \coloneqq 2 (f_{\A}(x),0,1) - f_{\mathrm{aux},\A}(y),
	\end{align}
then for $y\in R_\A$ it holds that $(f_{\A}(x),0,1) \in F_{\mathrm{TF}}(y)=\overline{\mathrm{conv}}\{f_{\mathrm{aux},\A}(y),f_{\mathrm{DF},\A}(y)\}$.
Similarly,	let the region $R_4$ be equipped with the following ODE
	\begin{align}\label{eq:dae_forming_ode2}
		y' &= f_{\mathrm{DF},\B}(y) \coloneqq 2 (f_{\B}(x),0,1) - f_{\mathrm{aux},\B}(y),
	\end{align}
then for $y\in R_\B$ it holds that $(f_{\B}(x),0,1) \in F_{\mathrm{TF}}(y)=\overline{\mathrm{conv}}\{f_{\mathrm{aux},\B}(y),f_{\mathrm{DF},\B}(y)\}$.
\end{proposition}
\textit{Proof.} We prove the assertion for Eq. \eqref{eq:dae_forming_ode1} and the second part follows similar lines. Note that for $y \in R_{\A} = \{ y \mid c(y)\coloneqq w = 0, \psi(x) < 1 \}$ we have that $\nabla c(y)^\top f_{\mathrm{aux},\A}(y)< 0$ and $\nabla c(y)^\top f_{\mathrm{DF},\A}(y)>0$. Hence, we have a sliding mode on $w =0$ with  $\frac{\mathrm{d} w}{\mathrm{d} \tau} =  0$ \cite{Filippov2013}.
From \eqref{eq:FilippovDI_with_multiplers} we have that $F_{\mathrm{TF}}(y) =\{ \theta_1 (2 (f_{\A}(x),0,1) - f_{\mathrm{aux},\A}(y)) + \theta_2 f_{\mathrm{aux},\A}(y) \mid \theta_1+\theta_2 =1 ,\theta_1,\theta_2 \geq 0\}$. From this relation and $w' =0$  we obtain that $\theta_1 - \theta_2 =0$. Thus we can solve for $\theta_1$ and $\theta_2$, i.e., $\theta_1 = \theta_2 = \frac{1}{2}$, which yields $(f_{\A}(x),0,1) \in F_{\mathrm{TF}}(y)$.
This completes the proof. \qed

Note that by construction the two sliding modes on $R_\A$ and $R_\B$ agree with the r.h.s. of Eq. \eqref{eq:ode_hysteresis_with_w} augmented by the dynamics of the clock state. Now we have defined vector fields in all regions of the time-freezing PSS which corresponds to the original system \eqref{eq:ode_hysteresis_with_w}. Another favorable property of the chosen auxiliary and DAE forming ODE is: since $w'(\tau)$ is bounded by $a>0$ it cannot make the sliding mode DAE arbitrarily stiff, especially if constraint drift happens.

\subsection{A Tutorial Example}
To illustrate the theoretical development we construct a time-freezing PSS for a thermostat system with hysteresis. The source code of the example is available in the repository of \texttt{NOSNOC} \cite{Nurkanovic2022c}.
 The system has a single state $x(\cdot)$ which models the temperature of a room which should stay inside the interval $x\in [18,20]$. As soon as the temperature drops below $x = 18$ the heater is switched on and when the temperature grows above $x = 20$ it is switched off.  The two modes of operation are given by $\dot{x} = f_{\A}(x) = -0.2x +5$ when the heater is on and $\dot{x} = f_{\B}(x) = -0.2x $ when the heater is off. One can see that for $\psi(x) = 0.5(x-18)$ we have a hybrid system that maches the finite automaton in Fig. \ref{fig:finite_automation}.
\begin{figure}[t]
	\centering
	{ \input{img/termostat_time_freezing_conf.tikz}}
	\vspace{-0.66cm}
	\caption{Trajectories of the time-freezing PSS for a thermostat example in numerical time $\tau$ (left plot) and physical time $t$ (right plot).}
	\vspace{-0.35cm}
	\label{fig:tutorail_example1}						
\end{figure}
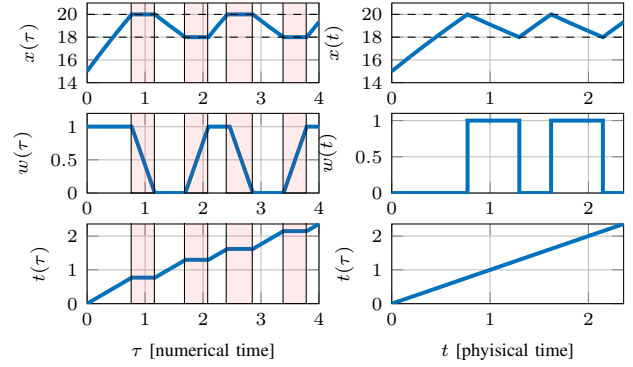
For a time-freezing PSS we define the regions $R_i$ via the Voronoi points as in the last section. The auxiliary ODE's r.h.s. according to Proposition \ref{prop:auxiliary_ode} read as $f_{\mathrm{aux},\A}(y) = (0,-\gamma(0.5(x-18)-1),0)$ and $f_{\mathrm{aux},\B}(y) = (0,\gamma(0.5(x-18),0)$ with $a = 1$. Similarly, the DAE-forming ODE r.h.s. according to Proposition \ref{prop:dae_forming_ode} read as $f_{\mathrm{DF},\A}(y) = (-0.4x +10,\gamma(0.5(x-18)-1),2)$ and $f_{\mathrm{aux},\B}(y) = (-0.4x,-\gamma(0.5(x-18)),2)$.

We simulate now the time-freezing PSS with a FESD Radau-IIA integrator of order 3 \cite{Nurkanovic2022} with $x(0) = 15$ and $w(0) = 0$.
The left plot in Fig. \ref{fig:tutorail_example1} illustrates the evolution of the time-freezing PSS in numerical time. The red shaded areas indicate the phases when the auxiliary ODE is active with $w \notin \{0,1\}$ while the time is frozen, cf. bottom left plot. In the middle left plot we can see that $w(\tau)$ is now a continuous function in numerical time. The right plot in Fig. \ref{fig:tutorail_example1} shows the differential state in physical time $t(\tau)$. Clearly, in the middle right plot $w(t(\tau))$ is now a discontinuous function, hence the state jumps are successfully recovered in physical time.


\section{Solution Equivalence}\label{sec:solution_relationship}
From the developments in the last section, the solution equivalence is nearly apparent. We formalize it in the next theorem.

\begin{theorem} \label{thm:solution_equivalence}
Regard the IVP corresponding to: 
(i) the Filippov DI of the time-freezing PSS equipped with the vector fields from Proposition \ref{prop:auxiliary_ode} and \ref{prop:dae_forming_ode} with a initial value $y(0) = (z_0,0)$ with $z_0 = (x_0,w_0)$ and $w_0\in \{0,1\}$, 
 on a time interval $[0,\tauf]$, 
(ii) the ODE with state jumps from Eq. \eqref{eq:ode_hysteresis_with_w} with $z(0)=z_0$ on a time interval $[0,\tf]  = [0,t(\tauf)]$. Suppose solutions exist to both IVP.
 
Then the solutions of the two IVPs $z(t;z_0)$ and $y(\tau;y_0)$ fulfill at any $\dfrac{dt}{d\tau}=t'(\tau)\neq 0$:
	\begin{align}\label{eq:solution_relation}
		z(t(\tau);z_0) &=\! My(t(\tau);y_0), 
\!		\text{ with }\!
		M \!\!=\! \begin{bmatrix}
			I_{n_x+1}\!\!\!  &\!\!\! \textbf{0}_{n_x+1,1} 
		\end{bmatrix}\!\!.
	\end{align}
\color{black}
\vspace{-0.5cm}
\end{theorem} 
\emph{Proof.}
Denote the solution of IVP (i) by $y_1(\tau;y_0)$ for $\tau \in (0,\hat{\tau})$ and for (ii) and $t(\tau) \in (0,t(\hat{\tau}))$ by $z_1(t(\tau);z_0)$. For a given $w(0) = \wone$ (or $\wtwo$) we have from Proposition \ref{prop:dae_forming_ode} that $y' = (f_\A(x),0,1)$ (or $y' = (f_\B(x),0,1)$).
Note that if there is no $\taus \in (0,\hat{\tau})$ for the IVP (i) such that an auxiliary ODE becomes active, then $t(\tau) = \int_{0}^{\tau} \dd \tau_1 = \tau$. Since $(f_\A(x),0) = M  (f_\A(x),0,1)$, $(f_\B(x),0) = M  (f_\B(x),0,1)$ and $z_0 = My_0$ by setting $\hat{\tau} = \tau_{\textrm{f}}$, it follows that \eqref{eq:solution_relation} holds.

Suppose now that we have a $\taus \in (0,\tauf)$ such that for $w(\taus) = \wtwo$ the auxiliary ODE $y' = f_{\mathrm{aux},\A}(y)$ becomes active (or similarly for $w(\taus) = \wone$, $y' = f_{\mathrm{aux},\B}(y)$ becomes active). From the first part of the proof we have that \eqref{eq:solution_relation} holds for $\tau \in(0,\tau_{\textrm{s}})$ and hence for all $t(\tau) \in(0,t_{\mathrm{s}}^-)$, where $t_{\mathrm{s}}^- = t(\taus)$. From Proposition \ref{prop:auxiliary_ode} we have that the solution satisfies $x(\taus)=x(\taur)$ and $w(\taur) = \wone$ (or $w(\taur) = \wtwo$) with $t'(\tau) = 0$ for $\tau \in [\taus,\taur]$. Hence, we have also  $t(\tau_{\textrm{r}}) = t_\mathrm{s}^+ = t(\tau_{\textrm{s}})$. Denote by $y_{\mathrm{s}} = (x(\taur),w(\taur),t(\taur))$.
Using this we have $y_1(\tau-\tau_{\textrm{r}},y_{\mathrm{s}}) = y(\tau,y_0)$ for $\tau \in(\tau_{\textrm{r}},\tilde{\tau})$ and denoting $z_{\mathrm{s}} = My_{\mathrm{s}}$ we see that $z_1(t(\tau)-t_s;z_{\mathrm{s}}) = x(t(\tau),z_0)$ for $t(\tau) \in (t_s^+,\tilde{\tau})$. Assume that a single activation of an auxiliary ODE takes place and set $\tilde{\tau} = \tauf$. 
Since the intervals $(t_s,t_{\textrm{f}})$ and $(\tau_{\textrm{r}},\tau_{\textrm{f}})$ have the same length and $z_{\mathrm{s}} = My_{\mathrm{s}}$ from the definitions of the corresponding IVP, we conclude that relation \eqref{eq:solution_relation} holds. If the auxiliary ODE becomes active multiple times we simply apply the same argument on the corresponding sub-intervals. This completes the proof. \qed

The last theorem opens the door to study the regarded hybrid system with hysteresis as a Filippov system and to apply their rich theory e.g., solution existence results \cite{Filippov2013}. 
From the practical side, we can use numerical methods for Filippov systems which allows us to avoid the use of integer variables.

\section{Numerical Example: Time Optimal Problem of a Car with Turbo Charger}\label{sec:numerical_example}

In this section we apply the theoretical developments in a numerical example of a time optimal control problem of a car with turbo from \cite{Avraam2000}. We consider a double-integrator car model equipped with a turbo accelerator which follows a hysteresis characteristic as in Fig. \ref{fig:hysteresis_char}. This makes the seemingly simple model severely nonlinear and nonsmooth. 

The car is described by its position $q(t)$, velocity $v(t)$ and turbo charger state $w(t) \in \{0,1\}$. The control variable is the car acceleration $u(t)$. The turbo accelerator is activated when the velocity exceeds $v \geq 15$ and is deactivated when it falls below $v \leq 10$. When it is on, it makes the nominal acceleration $u(t)$ three times greater. One can see that $\psi(x) = \frac{v-10}{5}$.
In summary, the state vector reads as $z =(q,v,w) \in \R^3$ with two modes of operation described by $f_{\A}(z) = (v,u,0)$ and $f_{\B}(z) = (v,3u,0)$. The acceleration is bounded by $|u| \leq \bar{u}$, $\bar{u} = 5$ and the velocity by $|v| \leq \bar{v}$, $\bar{v} = 25$. 

 
In the OCP we consider the time-freezing PSS associated to the car model on a numerical time interval $\tau \in [0,\tauf]$. The car should reach the goal $q(t(\tauf))\! = q_{\mathrm{f}} \!= 150$ with $v(t(\tauf))\!  = v_{\mathrm{f}}\! = 0$, whereby $z(0) = z_0 = \textbf{0}_{3,1}$.
The auxiliary and DAE-forming dynamics are chosen according to Propositions \ref{prop:auxiliary_ode} (with $a = 1$) and \ref{prop:dae_forming_ode}, respectively.
The OCP reads as:
\begin{subequations}\label{eq:ocp_time_freezing}
	\begin{align}
		\min_{y(\cdot),u(\cdot),s(\cdot)} \quad & t(\tauf) \\
		\textrm{s.t.} \quad & y(0) = (z_0,0),\\
		&{y}'(\tau)\! \in \! s(\tau)\! F_{\mathrm{TF}}(y(\tau),\!u(\tau)) ,\!\; \tau \in \![0,\tauf],\\
		& -\bar{u} \leq u(\tau) \leq \bar{u},\; \tau \in [0,\tauf],\\
		& \bar{s}^{-1} \leq s(\tau) \leq \bar{s},\;\tau \in [0,\tauf],\label{eq:speed_of_time}\\
		&-\bar{v} \leq v(\tau) \leq \bar{v},\;\tau \in [0,\tauf],\\
		&(q(\tauf),v(\tauf)) = (q_{\mathrm{f}},v_{\mathrm{f}}).
	\end{align}
\end{subequations}	

The objective consist of minimizing the final physical time. Since a time optimal control problem is considered, we introduce the scalar \textit{speed-of-time} control variable $s(\cdot)$ which introduces a time-transformation and enables to have a variable terminal physical time $T_{\mathrm{f}} = t(\tauf)$. It is bounded by \eqref{eq:speed_of_time} with $\bar{s} = 10$. \texttt{NOSNOC} ensures equidistant control grids in numerical and physical time $\tau$.

The OCP is discretized with a FESD Radau-IIA scheme of order 3 with $N=10$ control intervals and $N_{\mathrm{fe}} = 3$ additional integration steps on every control interval, with $\tauf = 5$. The controls are taken to be piecewise constant over the control intervals. The OCP discretization and MPCC homotopy is carried out via the open source tool \texttt{NOSNOC}, which has \texttt{IPOPT} \cite{Waechter2006} and \texttt{CasADi} \cite{Andersson2019} as a back-end.
 
Additionally, we compare our approach to the mixed integer formulation of \cite{Avraam2000}. We take the same control and state discretization as in \texttt{NOSNOC} which results $56$ binary variables. Switches in the integer formulation are allowed only at the control interval boundaries, as a switch detection formulation requires significantly more integer variables and introduces more nonlinearity.

 The problem is solved with the dedicated mixed integer nonlinear programming (MINLP) solver \texttt{Bonmin} \cite{Bonami2008}. Note that the only nonlinearity in the MINLP is due to the time-transformation for the optimal time $T_{\mathrm{f}}$. Therefore, in a second experiment we fix $T_{\mathrm{f}}$ and solve the resulting MILP with the commercial solver \texttt{Gurobi}. We make a bisection-type search in $T_{\mathrm{f}}$. The MILP with the smallest $T_{\mathrm{f}}$ that is still feasible, delivers the optimal time $T_{\mathrm{f}}$. In this experiment 22 MILP were solved for an accuracy of $10^{-6}$. To determine the solution quality, we additionally perform a high accuracy solution with the computed optimal controls and obtain $x_{\mathrm{sim}}(t)$. We compare the terminal constraint satisfactions: $E(T_{\mathrm{f}}) = \| x_{\mathrm{sim}}(T_{\mathrm{f}}) - (q_{\mathrm{f}},v_{\mathrm{f}}) \|_2$.  The source core for the simulation and the two MIOP approaches are provided in the \texttt{NOSNOC} repository \cite{Nurkanovic2022c}.

The results are summarized in Table \ref{tab:comparison_minlp}. All three approaches provide a similar objective value. \texttt{Gurobi} is the fastest solver, \texttt{NOSNOC} is only slightly slower and \texttt{Bonmin} is significantly slower. The smallest terminal error is achieved via \texttt{NOSNOC}. This is due to the underlying FESD discretization, cf. \cite{Nurkanovic2022}. \texttt{Gurobi} and \texttt{Bonmin} have the same discretization without switch detection and result in the same terminal error.
On the other hand, \texttt{Gurobi} provides the most robust approach, as \texttt{NOSNOC} (i.e., \texttt{IPOPT} as underlying NLP solver) fails to converge in some variations of the discretization. 
\color{black}
 \begin{figure}[t]
	\centering
	\vspace{0.25cm}
	{ \input{img/time_optimal_car1_conf.tikz}}
	\vspace{-0.6cm}
	\caption{Solutions of the OCP \eqref{eq:ocp_time_freezing} in physical time. The top left and right plots show the velocity $v(t)$ and optimal controls $u(t)$, respectively. The bottom left and right plots show the hysteresis state $w(t)$ and the solution trajectory in the $(v,w)$-plane, respectively.}
	\label{fig:time_optimal_car1}		
		\vspace{-0.3cm}				
\end{figure}
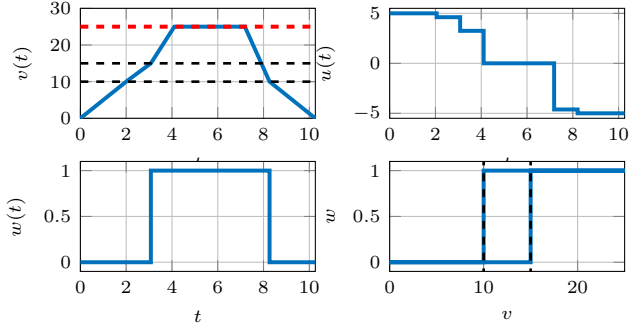


The results computed by \texttt{NOSNOC} is depicted in Fig. \ref{fig:time_optimal_car1}. One can see an intuitive behavior as the car uses the turbo accelerator as much as possible to reach the goal time optimally, with $T_{\mathrm{f}} = 10.26$. 



\section{Conclusion}\label{sec:conclusion}
In this paper we introduced a novel time-freezing reformulation for a class of hybrid systems with hysteresis. It transforms the systems with state jumps into PSS for which we leverage the recently developed FESD method which enables high accuracy optimal control by solving only smooth NLP. Thus, we can avoid use of computationally expensive mixed integer strategies in numerical optimal control and obtain quickly good and accurate nonsmooth solutions. In the theoretical part, constructive ways to find auxiliary and DAE-forming ODE are provided and solution equivalence is proven.  In future time-freezing for other types of finite automaton and hysteresis systems, as e.g., described in the introduction should be investigated as well.
\begin{table}[t]	
	\centering
	\vspace{0.4cm}
	\caption{Comparison of \texttt{NOSNOC} to mixed integer formulations.}
	\centering
	\begin{tabular}{|c|c|c|c|}	
		\hline
		Solver& {$T_{\mathrm{f}}$} & {CPU Time [s]} &{$E(T_{\mathrm{f}})$} \\
		\hline
		\texttt{NOSNOC} & \textbf{10.26} & 8.87 & \textbf{9.49e-02}  \\
		\hline
		\texttt{Gurobi} with bisection & 11.21  &\textbf{ 5.31} &7.88e+01  \\
		\hline
		\texttt{Bonmin} & 11.28 & 1481.58 & 7.88e+01 \\
		\hline
	\end{tabular}
	\label{tab:comparison_minlp}
\end{table}
\section*{Acknowledgment}
We thank Costas Pantelides from the Imperial College London and PSE Enterprise for inspiring discussions and providing the example in Section \ref{sec:numerical_example}.
\bibliographystyle{ieeetran}
\bibliography{bib/syscop1}
\end{document}

%% file: img/finite_automaton2.tikz
\pgfplotsset{compat=1.13}
\definecolor{mycolor1}{rgb}{0.00000,0.44700,0.74100}%
\definecolor{mycolor1}{rgb}{0.00000,0.0,0.0}%
\usetikzlibrary{automata, positioning, arrows}
\begin{tikzpicture}
\hspace{-1.2cm}
\node[state,text width=1.7cm, align = left] (1) 
{
$\dot{x} = f_{\mathrm{A}}(x)$ \\
$ w =0$ 
};
\node[state, right of=1, xshift=2.0cm,text width=1.7cm, align = left] (2) {
$\dot{x} = f_{\mathrm{B}}(x)$ \\
$ w = 1$ 
};
  \path[->]
(1) edge[bend left, above] node{\small$\psi(x)\geq 1$} (2) 
(2)	edge[bend left, below] node{\small$\psi(x)\leq 0$} (1);
\end{tikzpicture}

%% file: img/hysteresis_char.tikz
\pgfplotsset{compat=1.13}
\setlength{\fwidth}{5.5cm}
\setlength{\fheight}{1.9cm}
\definecolor{mycolor1}{rgb}{0.00000,0.44700,0.74100}%
\definecolor{mycolor1}{rgb}{0.00000,0.0,0.0}%
\begin{tikzpicture}

\begin{axis}[%
width=0.951\fwidth,
height=\fheight,
at={(0\fwidth,0\fheight)},
scale only axis,
xmin=-2.7,
xmax=2.3,
xtick={-1.5,0.75},
xticklabels={{$0$},{$1$}},
xlabel style={font=\color{white!15!black}},
xlabel={$\psi(x)$},
ymin=-0.1,
ymax=1.6,
ytick={0,0.2,0.4,0.6,0.8,1,1.2,1.4,1.6},
yticklabels={{},{},{},{},{},{},{},{},{}},
ytick={0,1.5},
yticklabels={{$0$},{$1$}},
ylabel style={font=\color{white!15!black}},
ylabel={$w$},
axis background/.style={fill=white},
legend style={legend cell align=left, align=left, draw=white!15!black},
xlabel style={font={\scriptsize}},ylabel style={font=\scriptsize},  ylabel shift={-0cm},ticklabel style={font=\scriptsize}
]
\addplot [color=mycolor1, line width=1.5pt]
  table[row sep=crcr]{%
-1.5	1.5\\
2.75	1.5\\
};

\addplot [color=mycolor1, line width=1.5pt]
  table[row sep=crcr]{%
-2.75	0\\
0.75	0\\
};

\addplot [color=mycolor1, draw=none, mark size=2.5pt, mark=*, mark options={solid, fill=mycolor1, mycolor1}]
  table[row sep=crcr]{%
0.75	0\\
};

\addplot [color=mycolor1, draw=none, mark size=2.5pt, mark=*, mark options={solid, fill=mycolor1, mycolor1}]
  table[row sep=crcr]{%
-1.5	1.5\\
};

\addplot [color=mycolor1, dashed, line width=1.5pt]
  table[row sep=crcr]{%
-1.5	0\\
-1.5	1.5\\
};

\addplot [color=mycolor1, dashed, line width=1.5pt]
  table[row sep=crcr]{%
0.75	0\\
0.75	1.5\\
};

\end{axis}
\tikzset{>=latex}
\draw[dashed,->=stealth,color=mycolor1,line width=1.0pt] (0.9,1.65) -- (0.9,0.65);
\draw[dashed,->=stealth,color=mycolor1,line width=1.0pt] (4.0,0.65) -- (4.0,1.65);

\begin{axis}[%
width=1.227\fwidth,
height=1.227\fheight,
at={(-0.16\fwidth,-0.135\fheight)},
scale only axis,
xmin=0,
xmax=1,
ymin=0,
ymax=1,
axis line style={draw=none},
ticks=none,
axis x line*=bottom,
axis y line*=left,
legend style={legend cell align=left, align=left, draw=white!15!black},
xlabel style={font={\scriptsize}},ylabel style={font=\scriptsize},  ylabel shift={-0cm},ticklabel style={font=\scriptsize}
]
\end{axis}
\end{tikzpicture}%

%% file: img/termostat_time_freezing_conf.tikz
\pgfplotsset{compat=1.13}
\setlength{\fwidth}{7.5cm}
\setlength{\fheight}{4cm}
\definecolor{mycolor1}{rgb}{0.00000,0.44700,0.74100}%
\begin{tikzpicture}

\begin{axis}[%
width=0.411\fwidth,
height=0.265\fheight,
at={(0\fwidth,0.735\fheight)},
scale only axis,
xmin=0,
xmax=4,
ymin=14,
ymax=21,
ylabel style={font=\color{white!15!black}},
ylabel={$x(\tau)$},
axis background/.style={fill=white},
xmajorgrids,
ymajorgrids,
legend style={legend cell align=left, align=left, draw=white!15!black},
xlabel style={font={\scriptsize}},ylabel style={font=\scriptsize},  ylabel shift={-0cm},ticklabel style={font=\scriptsize}
]
\addplot [color=mycolor1, line width=1.5pt]
  table[row sep=crcr]{%
0	15\\
0.0199999999999996	15.1397203750858\\
0.0399999999999991	15.27888298494\\
0.0599999999999987	15.4174900561758\\
0.0799999999999983	15.5555438065065\\
0.100000000000001	15.6930464447963\\
0.120000000000001	15.8300001710902\\
0.140000000000001	15.9664071766512\\
0.18	16.2375897467834\\
0.199999999999999	16.3723696504475\\
0.219999999999999	16.5066115112706\\
0.239999999999998	16.6403174771702\\
0.260000000000002	16.7734896874309\\
0.300000000000001	17.0382413553977\\
0.32	17.1698250493642\\
0.34	17.300883459779\\
0.359999999999999	17.4314186838436\\
0.379999999999999	17.5614328095116\\
0.399999999999999	17.6909279177255\\
0.420000000000002	17.8199060800276\\
0.440000000000001	17.9483693602132\\
0.460000000000001	18.0763198136726\\
0.48	18.2037594876588\\
0.5	18.3306904212479\\
0.52	18.4571146454099\\
0.539999999999999	18.583034183053\\
0.559999999999999	18.7084510491031\\
0.580000000000002	18.8333672506275\\
0.600000000000001	18.9577847870819\\
0.768615506136424	20.0000000028212\\
1.16861549792086	19.9999999913332\\
1.44	18.9433959418033\\
1.69350319664867	18.0000000011263\\
2.09350318987131	18.0000000145402\\
2.12	18.1691310497656\\
2.14	18.2962002205395\\
2.16	18.4227621298355\\
2.18	18.5488188029978\\
2.2	18.6743722570266\\
2.22	18.7994245011423\\
2.24	18.9239775368237\\
2.42	19.9994285556497\\
2.46	19.9999999480435\\
2.86	19.9999999330823\\
3.12932165579899	18.9512129347153\\
3.38	18.0245060457428\\
3.38903854457327	18.0000000038542\\
3.78903853469306	17.9999999966578\\
3.82	18.197541119827\\
3.84	18.3244968773413\\
3.86	18.4509458261335\\
3.88	18.5768899897668\\
3.9	18.7023313837141\\
3.92	18.8272720145601\\
3.94	18.9517138832875\\
4	19.3220667697353\\
};

\addplot[area legend, draw=black, fill=red, fill opacity=0.08]
table[row sep=crcr] {%
x	y\\
0.76	0\\
1.16	0\\
1.16	23.0000000028212\\
0.76	23.0000000028212\\
}--cycle;

\addplot[area legend, draw=black, fill=red, fill opacity=0.08]
table[row sep=crcr] {%
x	y\\
1.68	0\\
2.08	0\\
2.08	23.0000000028212\\
1.68	23.0000000028212\\
}--cycle;

\addplot[area legend, draw=black, fill=red, fill opacity=0.08]
table[row sep=crcr] {%
x	y\\
2.4	0\\
2.84931672107843	0\\
2.84931672107843	23.0000000028212\\
2.4	23.0000000028212\\
}--cycle;

\addplot[area legend, draw=black, fill=red, fill opacity=0.08]
table[row sep=crcr] {%
x	y\\
3.38	0\\
3.78	0\\
3.78	23.0000000028212\\
3.38	23.0000000028212\\
}--cycle;

\addplot [color=black, dashed]
  table[row sep=crcr]{%
0	18\\
0.0199999999999996	18\\
0.0399999999999991	18\\
0.0599999999999987	18\\
0.0799999999999983	18\\
0.100000000000001	18\\
0.120000000000001	18\\
0.140000000000001	18\\
0.18	18\\
0.199999999999999	18\\
0.219999999999999	18\\
0.239999999999998	18\\
0.260000000000002	18\\
0.300000000000001	18\\
0.32	18\\
0.34	18\\
0.359999999999999	18\\
0.379999999999999	18\\
0.399999999999999	18\\
0.420000000000002	18\\
0.440000000000001	18\\
0.460000000000001	18\\
0.48	18\\
0.5	18\\
0.52	18\\
0.539999999999999	18\\
0.559999999999999	18\\
0.580000000000002	18\\
0.600000000000001	18\\
2.12	18\\
2.14	18\\
2.16	18\\
2.18	18\\
2.2	18\\
2.22	18\\
2.24	18\\
2.46	18\\
3.82	18\\
3.84	18\\
3.86	18\\
3.88	18\\
3.9	18\\
3.92	18\\
3.94	18\\
4	18\\
};

\addplot [color=black, dashed]
  table[row sep=crcr]{%
0	20\\
0.0199999999999996	20\\
0.0399999999999991	20\\
0.0599999999999987	20\\
0.0799999999999983	20\\
0.100000000000001	20\\
0.120000000000001	20\\
0.140000000000001	20\\
0.18	20\\
0.199999999999999	20\\
0.219999999999999	20\\
0.239999999999998	20\\
0.260000000000002	20\\
0.300000000000001	20\\
0.32	20\\
0.34	20\\
0.359999999999999	20\\
0.379999999999999	20\\
0.399999999999999	20\\
0.420000000000002	20\\
0.440000000000001	20\\
0.460000000000001	20\\
0.48	20\\
0.5	20\\
0.52	20\\
0.539999999999999	20\\
0.559999999999999	20\\
0.580000000000002	20\\
0.600000000000001	20\\
2.12	20\\
2.14	20\\
2.16	20\\
2.18	20\\
2.2	20\\
2.22	20\\
2.24	20\\
2.46	20\\
3.82	20\\
3.84	20\\
3.86	20\\
3.88	20\\
3.9	20\\
3.92	20\\
3.94	20\\
4	20\\
};

\end{axis}

\begin{axis}[%
width=0.411\fwidth,
height=0.265\fheight,
at={(0\fwidth,0.368\fheight)},
scale only axis,
xmin=0,
xmax=4,
ymin=0,
ymax=1.2,
ylabel style={font=\color{white!15!black}},
ylabel={$w(\tau)$},
axis background/.style={fill=white},
xmajorgrids,
ymajorgrids,
legend style={legend cell align=left, align=left, draw=white!15!black},
xlabel style={font={\scriptsize}},ylabel style={font=\scriptsize},  ylabel shift={-0cm},ticklabel style={font=\scriptsize}
]
\addplot [color=mycolor1, line width=1.5pt]
  table[row sep=crcr]{%
0	1\\
0.0199999999999996	0.999999999996399\\
0.768615506136426	1.00000000104589\\
1.16861549792086	4.68782612728091e-09\\
1.69350319664867	3.07576186742153e-12\\
2.09350318987131	0.999999994879998\\
2.46	0.999999973467344\\
2.86	-1.32710393785374e-08\\
3.38903854457327	1.31604327435753e-09\\
3.78903853469306	0.999999995865208\\
4	1.00000000001366\\
};

\addplot[area legend, draw=black, fill=red, fill opacity=0.08]
table[row sep=crcr] {%
x	y\\
0.76	0\\
1.16	0\\
1.16	4.00000001327104\\
0.76	4.00000001327104\\
}--cycle;

\addplot[area legend, draw=black, fill=red, fill opacity=0.08]
table[row sep=crcr] {%
x	y\\
1.68	0\\
2.08	0\\
2.08	4.00000001327104\\
1.68	4.00000001327104\\
}--cycle;

\addplot[area legend, draw=black, fill=red, fill opacity=0.08]
table[row sep=crcr] {%
x	y\\
2.4	0\\
2.84931672107843	0\\
2.84931672107843	4.00000001327104\\
2.4	4.00000001327104\\
}--cycle;

\addplot[area legend, draw=black, fill=red, fill opacity=0.08]
table[row sep=crcr] {%
x	y\\
3.38	0\\
3.78	0\\
3.78	4.00000001327104\\
3.38	4.00000001327104\\
}--cycle;

\end{axis}

\begin{axis}[%
width=0.411\fwidth,
height=0.265\fheight,
at={(0\fwidth,0\fheight)},
scale only axis,
xmin=0,
xmax=4,
xlabel style={font=\color{white!15!black}},
xlabel={$\tau$ [numerical time]},
ymin=0,
ymax=2.35801385987606,
ylabel style={font=\color{white!15!black}},
ylabel={$t(\tau)$},
axis background/.style={fill=white},
xmajorgrids,
ymajorgrids,
legend style={legend cell align=left, align=left, draw=white!15!black},
xlabel style={font={\scriptsize}},ylabel style={font=\scriptsize},  ylabel shift={-0cm},ticklabel style={font=\scriptsize}
]
\addplot [color=mycolor1, line width=1.5pt]
  table[row sep=crcr]{%
0	0\\
0.0199999999999996	0.0200000000020459\\
0.04	0.0400000000034275\\
0.0599999999999996	0.0600000000053829\\
0.0800000000000001	0.0800000000074972\\
0.0999999999999996	0.100000000009967\\
0.12	0.120000000012775\\
0.14	0.14000000001599\\
0.18	0.180000000017022\\
0.2	0.200000000025537\\
0.22	0.220000000029829\\
0.24	0.240000000036604\\
0.26	0.260000000043825\\
0.3	0.300000000047968\\
0.32	0.320000000066115\\
0.34	0.34000000007789\\
0.36	0.360000000124811\\
0.38	0.380000000110726\\
0.4	0.400000000146421\\
0.42	0.420000000173101\\
0.44	0.440000000212767\\
0.46	0.460000000261669\\
0.48	0.480000000326651\\
0.5	0.500000000413968\\
0.52	0.520000000535324\\
0.54	0.540000000709878\\
0.56	0.560000000971827\\
0.58	0.580000001385526\\
0.600000000000001	0.600000002080523\\
0.768615506136426	0.770753845306196\\
1.16861549792086	0.770753838045258\\
1.7	1.29755676886256\\
2.09350318987131	1.29755676488955\\
2.12	1.32405357297831\\
2.14	1.34405357307999\\
2.16	1.36405357318305\\
2.18	1.38405357334325\\
2.2	1.40405357357574\\
2.22	1.42405357393937\\
2.24	1.44405357453865\\
2.42	1.62015412581191\\
2.46	1.62024935649532\\
2.86	1.6202493502716\\
3.38903854457327	2.14705240048961\\
3.78903853469306	2.14705239313396\\
3.82	2.17801385680486\\
3.84	2.1980138569134\\
3.86	2.21801385702474\\
3.88	2.23801385719892\\
3.9	2.25801385749472\\
3.92	2.2780138578365\\
3.94	2.29801385853551\\
4	2.35801385987606\\
};

\addplot[area legend, draw=black, fill=red, fill opacity=0.08]
table[row sep=crcr] {%
x	y\\
0.76	0\\
1.16	0\\
1.16	5.35801385987606\\
0.76	5.35801385987606\\
}--cycle;

\addplot[area legend, draw=black, fill=red, fill opacity=0.08]
table[row sep=crcr] {%
x	y\\
1.68	0\\
2.08	0\\
2.08	5.35801385987606\\
1.68	5.35801385987606\\
}--cycle;

\addplot[area legend, draw=black, fill=red, fill opacity=0.08]
table[row sep=crcr] {%
x	y\\
2.4	0\\
2.84931672107843	0\\
2.84931672107843	5.35801385987606\\
2.4	5.35801385987606\\
}--cycle;

\addplot[area legend, draw=black, fill=red, fill opacity=0.08]
table[row sep=crcr] {%
x	y\\
3.38	0\\
3.78	0\\
3.78	5.35801385987606\\
3.38	5.35801385987606\\
}--cycle;

\end{axis}

\begin{axis}[%
width=0.411\fwidth,
height=0.265\fheight,
at={(0.54\fwidth,0.735\fheight)},
scale only axis,
xmin=0,
xmax=2.35801385987606,
ymin=14,
ymax=21,
ylabel style={font=\color{white!15!black}},
ylabel={$x(t)$},
axis background/.style={fill=white},
xmajorgrids,
ymajorgrids,
legend style={legend cell align=left, align=left, draw=white!15!black},
xlabel style={font={\scriptsize}},ylabel style={font=\scriptsize},  ylabel shift={-0cm},ticklabel style={font=\scriptsize}
]
\addplot [color=mycolor1, line width=1.5pt]
  table[row sep=crcr]{%
0	15\\
0.0200000000020459	15.1397203750858\\
0.0400000000034275	15.27888298494\\
0.0600000000053811	15.4174900561758\\
0.0800000000074981	15.5555438065065\\
0.100000000009967	15.6930464447963\\
0.120000000012777	15.8300001710902\\
0.140000000015991	15.9664071766512\\
0.180000000017021	16.2375897467834\\
0.200000000025536	16.3723696504475\\
0.220000000029827	16.5066115112706\\
0.240000000036602	16.6403174771702\\
0.260000000043824	16.7734896874309\\
0.300000000047966	17.0382413553977\\
0.320000000066116	17.1698250493642\\
0.34000000007789	17.300883459779\\
0.36000000012481	17.4314186838436\\
0.380000000110726	17.5614328095116\\
0.40000000014642	17.6909279177255\\
0.420000000173101	17.8199060800276\\
0.440000000212766	17.9483693602132\\
0.460000000261669	18.0763198136726\\
0.480000000326651	18.2037594876588\\
0.500000000413969	18.3306904212479\\
0.520000000535326	18.4571146454099\\
0.540000000709878	18.583034183053\\
0.560000000971826	18.7084510491031\\
0.580000001385525	18.8333672506275\\
0.600000002080524	18.9577847870819\\
0.770753845306196	20.0000000028212\\
0.971582762197094	19.2126033326834\\
1.18213833825833	18.4203378276654\\
1.29755676486685	17.9999999953323\\
1.32405357297831	18.1691310497656\\
1.34405357307999	18.2962002205395\\
1.36405357318305	18.4227621298355\\
1.38405357334326	18.5488188029978\\
1.40405357357574	18.6743722570266\\
1.42405357393937	18.7994245011423\\
1.44405357453865	18.9239775368237\\
1.62024935676175	19.9999999486628\\
1.82024935558784	19.2157886785777\\
2.02893112046994	18.4302975057654\\
2.14705239268487	17.9999999929017\\
2.17801385680486	18.197541119827\\
2.1980138569134	18.3244968773413\\
2.21801385702474	18.4509458261335\\
2.23801385719892	18.5768899897668\\
2.25801385749472	18.7023313837141\\
2.2780138578365	18.8272720145601\\
2.29801385853551	18.9517138832875\\
2.35801385987606	19.3220667697353\\
};

\addplot [color=black, dashed]
  table[row sep=crcr]{%
0	18\\
0.0200000000020459	18\\
0.0400000000034275	18\\
0.0600000000053811	18\\
0.0800000000074981	18\\
0.100000000009967	18\\
0.120000000012777	18\\
0.140000000015991	18\\
0.180000000017021	18\\
0.200000000025536	18\\
0.220000000029827	18\\
0.240000000036602	18\\
0.260000000043824	18\\
0.300000000047966	18\\
0.320000000066116	18\\
0.34000000007789	18\\
0.36000000012481	18\\
0.380000000110726	18\\
0.40000000014642	18\\
0.420000000173101	18\\
0.440000000212766	18\\
0.460000000261669	18\\
0.480000000326651	18\\
0.500000000413969	18\\
0.520000000535326	18\\
0.540000000709878	18\\
0.560000000971826	18\\
0.580000001385525	18\\
0.600000002080524	18\\
0.770753845306196	18\\
1.29755676898971	18\\
1.32405357297831	18\\
1.34405357307999	18\\
1.36405357318305	18\\
1.38405357334326	18\\
1.40405357357574	18\\
1.42405357393937	18\\
1.44405357453865	18\\
1.62024935676175	18\\
2.14705240048961	18\\
2.17801385680486	18\\
2.1980138569134	18\\
2.21801385702474	18\\
2.23801385719892	18\\
2.25801385749472	18\\
2.2780138578365	18\\
2.29801385853551	18\\
2.35801385987606	18\\
};

\addplot [color=black, dashed]
  table[row sep=crcr]{%
0	20\\
0.0200000000020459	20\\
0.0400000000034275	20\\
0.0600000000053811	20\\
0.0800000000074981	20\\
0.100000000009967	20\\
0.120000000012777	20\\
0.140000000015991	20\\
0.180000000017021	20\\
0.200000000025536	20\\
0.220000000029827	20\\
0.240000000036602	20\\
0.260000000043824	20\\
0.300000000047966	20\\
0.320000000066116	20\\
0.34000000007789	20\\
0.36000000012481	20\\
0.380000000110726	20\\
0.40000000014642	20\\
0.420000000173101	20\\
0.440000000212766	20\\
0.460000000261669	20\\
0.480000000326651	20\\
0.500000000413969	20\\
0.520000000535326	20\\
0.540000000709878	20\\
0.560000000971826	20\\
0.580000001385525	20\\
0.600000002080524	20\\
0.770753845306196	20\\
1.29755676898971	20\\
1.32405357297831	20\\
1.34405357307999	20\\
1.36405357318305	20\\
1.38405357334326	20\\
1.40405357357574	20\\
1.42405357393937	20\\
1.44405357453865	20\\
1.62024935676175	20\\
2.14705240048961	20\\
2.17801385680486	20\\
2.1980138569134	20\\
2.21801385702474	20\\
2.23801385719892	20\\
2.25801385749472	20\\
2.2780138578365	20\\
2.29801385853551	20\\
2.35801385987606	20\\
};

\end{axis}

\begin{axis}[%
width=0.411\fwidth,
height=0.265\fheight,
at={(0.54\fwidth,0.368\fheight)},
scale only axis,
xmin=0,
xmax=2.35801385987606,
ymin=-1.04589004124023e-09,
ymax=1.1,
ylabel style={font=\color{white!15!black}},
ylabel={$w(t)$},
axis background/.style={fill=white},
xmajorgrids,
ymajorgrids,
legend style={legend cell align=left, align=left, draw=white!15!black},
xlabel style={font={\scriptsize}},ylabel style={font=\scriptsize},  ylabel shift={-0cm},ticklabel style={font=\scriptsize}
]
\addplot [color=mycolor1, line width=1.5pt]
  table[row sep=crcr]{%
0	0\\
0.0200000000020459	3.60067531346431e-12\\
0.040000000003428	3.42925687846218e-12\\
0.0600000000053824	4.67403893367191e-12\\
0.0800000000074976	5.45119505090952e-12\\
0.100000000009967	6.42419450969101e-12\\
0.120000000012775	7.35678185037614e-12\\
0.14000000001599	8.32267588180002e-12\\
0.180000000017022	-5.61328761250479e-13\\
0.200000000025537	1.48454581960777e-11\\
0.220000000029829	1.11226583499047e-11\\
0.240000000036604	1.36881617152085e-11\\
0.260000000043825	1.42352796217438e-11\\
0.300000000047968	-1.47437617670221e-13\\
0.320000000066115	2.30624408459335e-11\\
0.34000000007789	1.73261405223002e-11\\
0.36000000012481	1.011164485476e-10\\
0.380000000110726	-3.26672022765706e-12\\
0.400000000146421	3.19171356011339e-11\\
0.420000000173101	2.35336194975844e-11\\
0.440000000212767	2.99165137107593e-11\\
0.460000000261668	3.26294546937334e-11\\
0.480000000326651	3.82534004472745e-11\\
0.500000000413968	4.56195081710575e-11\\
0.520000000535325	5.67856872635275e-11\\
0.540000000709878	7.39377448155665e-11\\
0.560000000971827	1.01428643262125e-10\\
0.580000001385526	1.47269307859688e-10\\
0.600000002080523	2.26555219029478e-10\\
0.770753845306196	-1.04589004124023e-09\\
0.770753838045259	0.999999995312174\\
1.29755676898971	0.999999999996925\\
1.29755676488955	5.12000219998754e-09\\
1.32405357297832	-4.44089209850063e-16\\
1.34405357307999	5.52509149542857e-11\\
1.36405357318305	4.9388493295055e-11\\
1.38405357334325	6.97464308530016e-11\\
1.40405357357574	9.21480669546781e-11\\
1.42405357393937	1.32425181931239e-10\\
1.44405357453865	2.00332195277042e-10\\
1.62024935676175	3.48551751905291e-08\\
1.6202493502716	1.00000001327104\\
2.14705240048961	0.999999998683957\\
2.14705239313396	4.1347916202028e-09\\
2.17801385680486	-4.44089209850063e-16\\
2.1980138569134	5.75681724512833e-11\\
2.21801385702474	5.20858911556843e-11\\
2.23801385719892	7.42192973746114e-11\\
2.25801385749472	1.69794400761702e-10\\
2.2780138578365	1.2220624512338e-10\\
2.29801385853551	2.28312035943645e-10\\
2.35801385987606	-1.36584077381485e-11\\
};

\end{axis}

\begin{axis}[%
width=0.411\fwidth,
height=0.265\fheight,
at={(0.54\fwidth,0\fheight)},
scale only axis,
xmin=0,
xmax=2.35801385987606,
xlabel style={font=\color{white!15!black}},
xlabel={$t$ [phyisical time]},
ymin=0,
ymax=2.35801385987606,
ylabel style={font=\color{white!15!black}},
ylabel={$t(\tau)$},
axis background/.style={fill=white},
xmajorgrids,
ymajorgrids,
legend style={legend cell align=left, align=left, draw=white!15!black},
xlabel style={font={\scriptsize}},ylabel style={font=\scriptsize},  ylabel shift={-0cm},ticklabel style={font=\scriptsize}
]
\addplot [color=mycolor1, line width=1.5pt]
  table[row sep=crcr]{%
0	0\\
0.0200000000020459	0.0200000000020459\\
0.040000000003428	0.040000000003428\\
0.0600000000053824	0.0600000000053824\\
0.0800000000074976	0.0800000000074976\\
0.100000000009967	0.100000000009967\\
0.120000000012775	0.120000000012775\\
0.14000000001599	0.14000000001599\\
0.180000000017022	0.180000000017022\\
0.200000000025537	0.200000000025537\\
0.220000000029829	0.220000000029829\\
0.240000000036604	0.240000000036604\\
0.260000000043825	0.260000000043825\\
0.300000000047968	0.300000000047968\\
0.320000000066115	0.320000000066115\\
0.34000000007789	0.34000000007789\\
0.36000000012481	0.36000000012481\\
0.380000000110726	0.380000000110726\\
0.400000000146421	0.400000000146421\\
0.420000000173101	0.420000000173101\\
0.440000000212767	0.440000000212767\\
0.460000000261668	0.460000000261668\\
0.480000000326651	0.480000000326651\\
0.500000000413968	0.500000000413968\\
0.520000000535325	0.520000000535325\\
0.540000000709878	0.540000000709878\\
0.560000000971827	0.560000000971827\\
0.580000001385526	0.580000001385526\\
0.600000002080523	0.600000002080523\\
0.770753845306196	0.770753845306196\\
1.29755676898971	1.29755676898971\\
1.32405357297832	1.32405357297832\\
1.34405357307999	1.34405357307999\\
1.36405357318305	1.36405357318305\\
1.38405357334325	1.38405357334325\\
1.40405357357574	1.40405357357574\\
1.42405357393937	1.42405357393937\\
1.44405357453865	1.44405357453865\\
1.62024935676175	1.62024935676175\\
2.14705240048961	2.14705240048961\\
2.17801385680486	2.17801385680486\\
2.1980138569134	2.1980138569134\\
2.21801385702474	2.21801385702474\\
2.23801385719892	2.23801385719892\\
2.25801385749472	2.25801385749472\\
2.2780138578365	2.2780138578365\\
2.29801385853551	2.29801385853551\\
2.35801385987606	2.35801385987606\\
};

\end{axis}

\begin{axis}[%
width=1.227\fwidth,
height=1.227\fheight,
at={(-0.16\fwidth,-0.135\fheight)},
scale only axis,
xmin=0,
xmax=1,
ymin=0,
ymax=1,
axis line style={draw=none},
ticks=none,
axis x line*=bottom,
axis y line*=left,
legend style={legend cell align=left, align=left, draw=white!15!black},
xlabel style={font={\scriptsize}},ylabel style={font=\scriptsize},  ylabel shift={-0cm},ticklabel style={font=\scriptsize}
]
\end{axis}
\end{tikzpicture}%

%% file: img/time_optimal_car1_conf.tikz
\pgfplotsset{compat=1.13}
\setlength{\fwidth}{7.6cm}
\setlength{\fheight}{3.5cm}

\definecolor{mycolor1}{rgb}{0.00000,0.44700,0.74100}%
\begin{tikzpicture}

\begin{axis}[%
width=0.411\fwidth,
height=0.418\fheight,
at={(0\fwidth,0.582\fheight)},
scale only axis,
xmin=0,
xmax=10.2628974060198,
xlabel style={font=\color{white!15!black}},
xlabel={$t$},
ymin=0,
ymax=30,
ylabel style={font=\color{white!15!black}},
ylabel={$v(t)$},
axis background/.style={fill=white},
xmajorgrids,
ymajorgrids,
legend style={legend cell align=left, align=left, draw=white!15!black},
xlabel style={font={\scriptsize}},ylabel style={font=\scriptsize},  ylabel shift={-0cm},ticklabel style={font=\scriptsize}
]
\addplot [color=mycolor1, line width=1.5pt]
  table[row sep=crcr]{%
0	0\\
2.05257948120397	10.2628974060198\\
3.07886922180595	15\\
4.10515896240794	25\\
7.18402818421389	25\\
8.26289740601984	10\\
10.2628974060198	0\\
};

\addplot [color=black, dashed, line width=1.0pt]
  table[row sep=crcr]{%
0	10\\
3.07886922180595	10\\
8.26289740601984	10\\
10.2628974060198	10\\
};

\addplot [color=black, dashed, line width=1.0pt]
  table[row sep=crcr]{%
0	15\\
3.07886922180595	15\\
8.26289740601984	15\\
10.2628974060198	15\\
};

\addplot [color=red, dashed, line width=1.5pt]
  table[row sep=crcr]{%
0	25\\
3.07886922180595	25\\
8.26289740601984	25\\
10.2628974060198	25\\
};

\end{axis}

\begin{axis}[%
width=0.411\fwidth,
height=0.418\fheight,
at={(0.541\fwidth,0.582\fheight)},
scale only axis,
xmin=0,
xmax=10.2628974060198,
xlabel style={font=\color{white!15!black}},
xlabel={$t$},
ymin=-5.5,
ymax=5.5,
ylabel style={font=\color{white!15!black}},
ylabel={$u(t)$},
axis background/.style={fill=white},
xmajorgrids,
ymajorgrids,
legend style={legend cell align=left, align=left, draw=white!15!black},
xlabel style={font={\scriptsize}},ylabel style={font=\scriptsize},  ylabel shift={-0cm},ticklabel style={font=\scriptsize}
]
\addplot[const plot, color=mycolor1, line width=1.5pt] table[row sep=crcr] {%
0	5\\
1.02628974060198	5\\
2.05257948120397	4.61575557717409\\
3.07886922180595	3.24794568381646\\
4.10515896240794	-3.5527136788005e-15\\
5.13144870300992	-1.77635683940025e-15\\
6.1577384436119	0\\
7.18402818421389	-4.6157555771741\\
8.21031792481587	-5\\
9.23660766541786	-5\\
10.2628974060198	-5\\
};

\end{axis}

\begin{axis}[%
width=0.411\fwidth,
height=0.418\fheight,
at={(0\fwidth,0\fheight)},
scale only axis,
xmin=0,
xmax=10.2628974060198,
xlabel style={font=\color{white!15!black}},
xlabel={$t$},
ymin=-0.1,
ymax=1.1,
ylabel style={font=\color{white!15!black}},
ylabel={$w(t)$},
axis background/.style={fill=white},
xmajorgrids,
ymajorgrids,
legend style={legend cell align=left, align=left, draw=white!15!black},
xlabel style={font={\scriptsize}},ylabel style={font=\scriptsize},  ylabel shift={-0cm},ticklabel style={font=\scriptsize}
]
\addplot [color=mycolor1, line width=1.5pt]
  table[row sep=crcr]{%
0	0\\
3.07886922180595	-0\\
3.07886922180595	1\\
8.26289740601984	1\\
8.26289740601984	-0\\
10.2628974060198	-0\\
};

\end{axis}

\begin{axis}[%
width=0.411\fwidth,
height=0.418\fheight,
at={(0.541\fwidth,0\fheight)},
scale only axis,
xmin=0,
xmax=25,
xlabel style={font=\color{white!15!black}},
xlabel={$v$},
ymin=-0.1,
ymax=1.1,
ylabel style={font=\color{white!15!black}},
ylabel={$w$},
axis background/.style={fill=white},
xmajorgrids,
ymajorgrids,
legend style={legend cell align=left, align=left, draw=white!15!black},
xlabel style={font={\scriptsize}},ylabel style={font=\scriptsize},  ylabel shift={-0cm},ticklabel style={font=\scriptsize}
]
\addplot [color=mycolor1, line width=1.5pt]
  table[row sep=crcr]{%
0	0\\
15	-0\\
15	1\\
25	1\\
10	1\\
10	-0\\
0	-0\\
};

\addplot [color=black, dashed, line width=1.0pt]
  table[row sep=crcr]{%
10	-0.5\\
10	1.5\\
};

\addplot [color=black, dashed, line width=1.0pt]
  table[row sep=crcr]{%
15	-0.5\\
15	1.5\\
};

\end{axis}

\begin{axis}[%
width=1.228\fwidth,
height=1.228\fheight,
at={(-0.16\fwidth,-0.136\fheight)},
scale only axis,
xmin=0,
xmax=1,
ymin=0,
ymax=1,
axis line style={draw=none},
ticks=none,
axis x line*=bottom,
axis y line*=left,
legend style={legend cell align=left, align=left, draw=white!15!black},
xlabel style={font={\scriptsize}},ylabel style={font=\scriptsize},  ylabel shift={-0cm},ticklabel style={font=\scriptsize}
]
\end{axis}
\end{tikzpicture}%